\documentclass[12pt]{article}
\usepackage{amssymb}
\usepackage{amsmath}
\usepackage{amsfonts}
\usepackage{cite}

\begin{document}

\author{I.P. Smirnov \\
Institute of Applied Physics RAS, \\
46 Ul'yanova Street, Nizhny Novgorod, Russia}
\title{Optimal control of a linear system with multiplicative noise at
control parameter}
\date{}
\maketitle

\begin{abstract}
We investigate a control process described by a linear system of ordinary
differential equations with a noise of special type acting to the control
parameter. As the cost functional the probability of the final state vector
to enter to a given set in the phase space is considered. Necessary
conditions of optimality (of the Pontryagin maximum principle form) and
existence theorems are developed. The initial control problem was trasformed
to an auxiliary deterministic problem, the differentiability of the
auxiliary functional was discussed.
\end{abstract}

\begin{center}
\end{center}

\section{Introduction}

\label{intro}

Consider the optimal control problem 
\begin{gather}
\left\{ 
\begin{array}{l}
\frac{d}{dt}\mathbf{x}=A\left( t\right) \mathbf{x+}\xi \left( t\right)
B\left( t\right) \mathbf{u}\left( t\right) \mathbf{+f}\left( t\right) ,\
0\leq t\leq 1, \\ 
\mathbf{x}\left( 0\right) =\mathbf{x}_{0},%
\end{array}%
\right.   \label{1} \\
\varphi \left( \mathbf{u}\left( \cdot \right) \right) =\mathbf{P}\left( 
\mathbf{x}_{u}\left( 1\right) \in \mathbf{Q}\right) \rightarrow \max_{%
\mathbf{u}\left( \cdot \right) \in \mathfrak{U}},  \label{1a}
\end{gather}%
where $\mathbf{x}$ is the system state vector and $\mathbf{u}$ is the
control one, $\mathbf{Q}$ is a given nonrandom goal set in the phase space, $%
\xi \left( t\right) $ is a scalar random process of the given type.

The necessary conditions for optimality in the form of Pontryagin maximum
principle and existence theorems are obtained in this paper keeping on our
recent investigations \cite{Smirn1,Smirn2,Smirn3}. In sec. \ref{putt} we
state the problem, in sec. \ref{reduc} reduce it to the equivalent one of
deterministic type. Existence of optimal control is proved in sec. \ref%
{exist}. Conditions of differentiability of the functional, necessary and
sufficient conditions for optimality are obtained in sec. \ref{dif}, \ref%
{nec}, respectively.

\section{Statement of the problem}

\label{putt}

Denote by $\left\langle \cdot,\cdot\right\rangle $ the inner product and by $%
\left\vert \cdot\right\vert $ the Euclidean norm of a vector (matrix) in $%
\mathbb{R}^{n}$. Let $(\Omega,\mathfrak{F},\mathbf{P)}$\textbf{\ }be a
complete probability space. Forming the system (\ref{1}) deterministic
matrices $A\left( t\right) $, $B\left( t\right) $, $f\left( t\right) $ are
of order $n\times n$, $n\times m$, $n\times1$, respectively; they have
summable with respect to Lebesque measure components; the ones of matrix $%
B\left( t\right) $ are quadratically summable; $\mathbf{x}_{0}$ is a
nonrandom vector, $Q$ is a nonrandom Borel's set in $\mathbb{R}^{n}$.

Let $\left\{ \xi_{1}\left( \omega\right) ,\ldots,\xi_{N}\left( \omega\right)
\right\} $ be a given set of random values on $(\Omega ,\mathfrak{F},\mathbf{%
P)}$, 
\begin{equation*}
\left\{ 0=t_{0}<t_{1}<\ldots<t_{N}=1\right\} 
\end{equation*}
a$\ $finite subset of $\left[ 0,1\right] $. Now we introduce the process $%
\xi\left( t,\omega\right) $: 
\begin{equation*}
\xi\left( t,\omega\right) =\xi_{i}\left( \omega\right) \ ,\ t\in\left[
t_{i-1},t_{i}\right) ,\ i=\overline{1,N}. 
\end{equation*}

The set $\mathfrak{U}$ of admissible controls consists of all measurable
deterministic functions $\mathbf{u}\left( t\right) ,\ t\in \left[ 0,1\right] 
$, taking their values in a fixed set $\mathbf{V}\subset \mathbb{R}^{n}$.
For every admissible control\ $\mathbf{u}\left( \cdot \right) \in L_{2}^{m}%
\left[ 0,1\right] $ there exist a unique (in the class of process with
absolutely continuous trajectories) solution of the Cauchy's problem (\ref{1}%
)\cite[p. 35]{Gihman}; the corresponding value of the goal functional $%
\varphi \left( \mathbf{u}\left( \cdot \right) \right) \in \left[ 0,1\right] $%
. Optimal control problem (\ref{1})-(\ref{1a}) can be posed as the problem
of choosing of an admissible control $\mathbf{u}^{0}\left( \cdot \right) \in 
\mathfrak{U}$ provided maximal value to the functional (\ref{1a}):%
\begin{equation*}
\varphi \left( \mathbf{u}^{0}\left( \cdot \right) \right) =\max_{\mathbf{u}%
\left( \cdot \right) \in \mathfrak{U}}\varphi \left( \mathbf{u}\left( \cdot
\right) \right) .
\end{equation*}

\section{Auxiliary deterministic problem}

\label{reduc}

Using Cauchy's formula for the solution of a linear differential equation
system one can write%
\begin{equation}
\begin{array}{c}
\mathbf{x}_{u}\left( 1\right) =\Phi\left( 1\right) \left( \mathbf{x}%
_{0}+\int\limits_{0}^{1}\Phi^{-1}\left( s\right) \left( \xi\left( s\right)
B\left( s\right) \mathbf{u}\left( s\right) +\mathbf{f}\left( s\right)
\right) ds\right) = \\ 
=\mathbf{\hat{x}}_{0}+\Phi\left( 1\right)
\sum\limits_{I=1}^{N}\xi_{i}\int\limits_{t_{i-1}}^{t_{i}}\Phi^{-1}\left(
s\right) B\left( s\right) \mathbf{u}\left( s\right) ds=\mathbf{\hat{x}}%
_{0}+\sum\limits_{I=1}^{N}\xi_{i}\mathbf{z}_{u}^{\left( i\right) },%
\end{array}
\label{2}
\end{equation}
where 
\begin{align}
\mathbf{\hat{x}}_{0} & \equiv\Phi\left( 1\right) \left( \mathbf{x}%
_{0}+\int\limits_{0}^{1}\Phi^{-1}\left( s\right) \mathbf{f}\left( s\right)
ds\right) ,  \notag \\
\mathbf{z}_{u}^{\left( i\right) } & \equiv\Phi\left( 1\right)
\int\limits_{t_{i-1}}^{t_{i}}\Phi^{-1}\left( s\right) B\left( s\right) 
\mathbf{u}\left( s\right) ds,   \label{zu}
\end{align}
$\Phi\left( t\right) $ is a fundamental matrix for the system (\ref{1}) ($%
\Phi\left( 0\right) \ $is the unitary matrix). Note that

\begin{equation}
\mathbf{z}_{u}^{\left( i\right) }=\Psi_{i}\mathbf{y}\left( t_{i}\right)
-\Psi_{i-1}\mathbf{y}\left( t_{i-1}\right) ,\ \Psi_{i}\equiv\Phi\left(
1\right) \Phi^{-1}\left( t_{i}\right) ,i=\overline{1,N},\   \label{3}
\end{equation}
where $\mathbf{y}_{u}\left( t\right) $ is the solution of the Cauchy's
problem%
\begin{equation*}
\left\{ 
\begin{array}{l}
\frac{d}{dt}\mathbf{y=}A\left( t\right) \mathbf{y+}B\left( t\right) \mathbf{%
u,\ }t\in\left( 0,1\right) , \\ 
\mathbf{y}\left( 0\right) =0.%
\end{array}
\right. 
\end{equation*}
Taking into account \ref{2}) and (\ref{3}) and introducing the functions of $%
N$ vector variables%
\begin{equation}
\begin{array}{c}
g\left( \mathbf{z}_{1},\ldots,\mathbf{z}_{N}\right) =\mathbf{P}\left(
\sum\limits_{I=1}^{N}\xi_{i}\mathbf{z}_{i}\in\mathbf{\hat{Q}}\right) , \\ 
G\left( \mathbf{y}_{1},\ldots,\mathbf{y}_{N}\right) =\mathbf{P}\left(
\sum\limits_{I=1}^{N}\chi_{i}\Psi_{i}\mathbf{y}_{i}\in\mathbf{\hat{Q}}%
\right) ,%
\end{array}
\label{4}
\end{equation}
where $\mathbf{\hat{Q}\equiv Q-x}_{0}$, $\chi_{i}\equiv\xi_{i}-\xi_{i+1}$, $%
i=\overline{1,N-1}$, $\chi_{N}\equiv\xi_{N}$, we get the following
representation for the functional (\ref{1a}):%
\begin{equation}
\begin{array}{c}
\varphi\left( \mathbf{u}\left( \cdot\right) \right) =\mathbf{P}\left( 
\mathbf{x}_{u}\left( 1\right) \in\mathbf{Q}\right) = \\ 
=\mathbf{P}\left( \mathbf{\hat{x}}_{0}+\sum\limits_{I=1}^{N}\xi_{i}\mathbf{z}%
_{u}^{\left( i\right) }\in\mathbf{Q}\right) = \\ 
=g\left( \mathbf{z}_{u}^{\left( 1\right) },\ldots,\mathbf{z}_{u}^{\left(
N\right) }\right) =G\left( \mathbf{y}_{u}\left( t_{1}\right) ,\ldots,\mathbf{%
y}_{u}\left( t_{N}\right) \right) .%
\end{array}
\label{5}
\end{equation}

It is clear now that the stochastic control problem (\ref{1}) is equivalent
to the next deterministic control problem: 
\begin{gather}
\left\{ 
\begin{array}{l}
\frac{d}{dt}\mathbf{y=}A\left( t\right) \mathbf{y+}B\left( t\right) \mathbf{%
u,\ }t\in\left( 0,1\right) , \\ 
\mathbf{y}\left( 0\right) =0,%
\end{array}
\right.  \label{detprob} \\
\varphi\left( \mathbf{u}\left( \cdot\right) \right) =G\left( \mathbf{y}%
_{u}\left( t_{1}\right) ,\ldots,\mathbf{y}_{u}\left( t_{N}\right) \right)
\rightarrow\max_{\mathbf{u}\left( \cdot\right) \in\mathfrak{U}}.  \notag
\end{gather}

Note that the form of the last is unusual: the cost functional is of
terminal type and it depends on control process states, taken at several
moments of the time of control.

\section{Existence theorem}

\label{exist}

Let $\left\{ \mathbf{y}_{1},\ldots,\mathbf{y}_{N}\right\} $ be an arbitrary
set of vectors in $\mathbb{R}^{n}$. Construct a new vector $\mathbf{Y}$ in $%
\mathbb{R}^{M}$, $M\equiv n\cdot N$ by the following rule: $\mathbf{Y=}%
\left\{ \mathbf{y}_{1},\ldots,\mathbf{y}_{N}\right\} =\left\{
y_{1}^{1},\ldots,y_{1}^{n},y_{2}^{1},\ldots,y_{N}^{n}\right\} $. Consider
the following functions 
\begin{equation*}
\begin{array}{l}
S\left( \omega,\mathbf{Y}\right) =\sum\limits_{I=1}^{N}\chi_{i}\left(
\omega\right) \Psi_{i}\mathbf{y}_{i}, \\ 
h\left( \mathbf{Y}\right) =\mathbf{P}\left( \omega:S\left( \omega ,\mathbf{Y}%
\right) \in\mathbf{\hat{Q}}\right) .%
\end{array}
\end{equation*}

\textbf{Lemma 1.} \textit{If the set }$\mathbf{\hat{Q}}$ \textit{is closed
in }$\mathbb{R}^{n}$, \textit{then }$h\left( \mathbf{\cdot}\right) $ \textit{%
is upper semicontinues function}.

\textbf{Prove.} Let $\mathbf{Y}^{\ast}\in\mathbb{R}^{M}$, $\mathbf{Y}%
_{k}\rightarrow\mathbf{Y}^{\ast}\ $as $k\rightarrow\infty$. Consider random
events

\begin{equation*}
\begin{array}{l}
A_{k}=\left\{ \omega:S\left( \omega,\mathbf{Y}_{k}\right) \in \mathbf{\hat{Q}%
}\right\} , \\ 
A^{\ast}=\left\{ \omega:S\left( \omega,\mathbf{Y}^{\ast}\right) \in\mathbf{%
\hat{Q}}\right\} , \\ 
B=~\underset{k\rightarrow\infty}{\lim\sup}A_{k}=\underset{k\geq1}{\cap }%
\underset{j\geq k}{\cup}A_{j}.%
\end{array}
\end{equation*}

If $\omega \in B$, then there exist a sequence $k_{m}\rightarrow \infty $
such that $\omega \in A_{k_{m}}$, in other words 
\begin{equation*}
S\left( \omega ,\mathbf{Y}_{k_{m}}\right) \in \mathbf{\hat{Q}}
\end{equation*}%
for evey $m$. Proceeding in this statement to the limit as $m\rightarrow
\infty $ and taking into account that the function $S\left( \omega ,\cdot
\right) $ is continuous and the set $\mathbf{\hat{Q}}$ is closed we obtain $%
S\left( \omega ,\mathbf{Y}^{\ast }\right) \in \mathbf{\hat{Q}}$, so $\omega
\in A^{\ast }$. Hence $B\subset A^{\ast }$ and $\ \mathbf{P}\left( B\right)
\leq \mathbf{P}\left( A\right) $. On the other hand, due to measure
properties \cite[p.34]{Neve} 
\begin{equation*}
\mathbf{P}\left( B\right) =\mathbf{P}\left( \underset{k\rightarrow \infty }{%
\lim \sup }~A_{k}\right) \geq ~\underset{k\rightarrow \infty }{\lim \sup }~%
\mathbf{P}\left( A_{k}\right) .
\end{equation*}

Finally we have%
\begin{equation*}
\begin{array}{c}
h\left( \mathbf{Y}^{\ast}\right) =\mathbf{P}\left( A^{\ast}\right) \geq%
\mathbf{P}\left( B\right) \geq~\underset{k\rightarrow\infty}{\lim\sup }~%
\mathbf{P}\left( A_{k}\right) = \\ 
=\underset{k\rightarrow\infty}{~\lim\sup~}h\left( \mathbf{Y}_{k}\right) .%
\end{array}
\end{equation*}
This proves the lemma.

\textbf{Theorem 1}\footnote{%
This theorem was proved by A.Yu. Zorin}\textbf{.}\textit{\ Suppose that }$%
\mathfrak{U}$\textit{\ is a weakly compact in }$L_{2}^{m}\left[ 0,1\right] $%
\textit{\ and the goal set }$Q$\textit{\ is closed in }$R^{n}$\textit{. Then
the optimal problem (\ref{1})-(\ref{1a})\ has a solution.}

\textbf{Prove.} Let $\left\{ \mathbf{u}_{l}\left( \cdot\right) \right\} $ be
a maximizing, weakly converging to $\mathbf{u}^{0}\left( \cdot\right) \in%
\mathfrak{U}$ sequence of controls: 
\begin{equation*}
\varphi\left( \mathbf{u}_{l}\left( \cdot\right) \right) \rightarrow
\sup\varphi\left( \mathbf{u}\left( \cdot\right) \right) \leq1. 
\end{equation*}
It is easy to prove that $\mathbf{y}_{u_{l}}\left( t_{k}\right) \rightarrow%
\mathbf{y}_{u^{0}}\left( t_{k}\right) $, $k=\overline{1,N}$. According to
lemma 1 we have%
\begin{equation*}
\begin{array}{c}
\varphi\left( \mathbf{u}^{0}\left( \cdot\right) \right) =G\left( \mathbf{y}%
_{u^{0}}\left( t_{1}\right) ,\ldots,\mathbf{y}_{u^{0}}\left( t_{N}\right)
\right) \geq \\ 
\geq~\underset{l\rightarrow\infty}{\lim\sup}~G\left( \mathbf{y}%
_{u_{l}}\left( t_{1}\right) ,\ldots,\mathbf{y}_{u_{l}}\left( t_{N}\right)
\right) = \\ 
=\underset{l\rightarrow\infty}{~\lim}\varphi\left( \mathbf{u}_{l}\left(
\cdot\right) \right) =\underset{\mathbf{u}\left( \cdot\right) \in\mathfrak{U}%
}{\sup}\varphi\left( \mathbf{u}\left( \cdot\right) \right) .%
\end{array}
\end{equation*}
This proves the theorem.

Another approaches to the existence of the solutions of such kind problems
were developed in \cite{Smirn1,Smirn2,Smirn3}.

\section{Conditions of differentiability}

\label{dif}

Differentiability of the functional with respect to the phase coordinates
plays important role in necessary and sufficient conditions of optimality.
Is this section we study the differentiability of functions $g,G$.

Restrict our consideration to the case of absolutely continuous with respect
to Lebesque measure random vector $\xi\equiv\left\{ \xi_{1},\ldots,\xi
_{N}\right\} $. Let $f$ be the probability density function of $\xi$, $%
\mathbf{Z\equiv}\left\{ \mathbf{z}_{1}^{\ast},\ldots,\mathbf{z}_{N}^{\ast
}\right\} $ be a given set of vectors in $\mathbb{R}^{n}$, 
\begin{equation*}
O\left( \mathbf{Z}\right) \equiv\left\{ \mathbf{r\in}\mathbb{R}^{N}:r_{1}%
\mathbf{z}_{1}^{\ast}+\ldots+r_{N}\mathbf{z}_{N}^{\ast}\in \mathbf{Q}%
\right\} . 
\end{equation*}
If $N>n$ then the set $O\left( \mathbf{Z}\right) $ is unbounded in $\mathbb{R%
}^{N}$ even $\mathbf{Q}$ is bounded in $\mathbb{R}^{n}$. Let $\pi_{j}$ be
the projection of $O\left( \mathbf{Z}\right) $ to the hyperplane $r_{j}=0$.
It is easily shown that for any convex $\mathbf{Q}$ the line passing through
an arbitrary point $\mathbf{p}_{j}\in\pi_{j}$ in parallel to the axe $r_{j}$
intersects $O\left( \mathbf{Z}\right) $ over interval $r_{j}\in\left(
r_{j}^{\left( 1\right) },r_{j}^{\left( 2\right) }\right) $; the functions $%
r_{j}=r_{j}^{\left( 1,2\right) }\left( \mathbf{p}_{j}\right) $ are bounded
for a bounded $\mathbf{Q}$ and $\mathbf{z}_{j}^{\ast}\neq0$.

For a fixed $k\in\overline{1,N}$ consider the function%
\begin{equation*}
h_{k}\left( \mathbf{z}_{k}\right) =g\left( \mathbf{z}_{1}^{\ast},\ldots,%
\mathbf{z}_{k}^{\ast},\mathbf{z}_{k},\mathbf{z}_{k+1}^{\ast},\ldots,\mathbf{z%
}_{N}^{\ast}\right) ,\ \mathbf{z}_{k}\in\mathbb{R}^{n}. 
\end{equation*}

\textbf{Lemma 2.} \textit{Let }$j,k\in\overline{1,N}$\textit{, vector }$%
z_{j}^{\ast}\neq0$\textit{, the probability density }$f\left(
r_{1},\ldots,r_{N}\right) $\textit{\ is differentiable with respect to }$%
r_{j}$\textit{\ and the function }$\Phi_{kj}\left( \mathbf{r}\right)
\equiv\partial\left( r_{k}f\right) /\partial r_{j}$\textit{\ is summable in }%
$\mathbb{R}^{N}$\textit{. Then the function }$h_{k}\left( \cdot\right) $%
\textit{\ is differentiable at the point }$\mathbf{z}_{k}^{\ast}$\textit{\
in the direction of vector }$\mathbf{z}_{j}^{\ast}$\textit{. The
corresponding directional derivative}

\begin{equation}
\frac{\partial h_{k}}{\partial\mathbf{z}_{j}^{\ast}}=\frac{d}{d\varepsilon }%
h_{k}\left( \mathbf{z}_{k}^{\ast}+\varepsilon\mathbf{z}_{j}^{\ast}\right)
\mid_{\varepsilon=0}=-\int\limits_{O\left( \mathbf{Z}\right) }\Phi
_{kj}\left( \mathbf{r}\right) d\mathbf{r.}   \label{6}
\end{equation}

\textbf{Prove.} We have%
\begin{equation}
h_{k}\left( \mathbf{z}_{k}^{\ast}+\varepsilon\mathbf{z}_{j}^{\ast}\right) =%
\mathbf{P}\left\{ \xi_{1}\mathbf{z}_{1}^{\ast}+\ldots+\xi_{k}\mathbf{z}%
_{k}^{\ast}+\ldots+\xi_{N}\mathbf{z}_{N}^{\ast}\in\mathbf{Q}\right\}
=\int\limits_{O\left( \mathbf{Z}\right) }f_{\eta}\left( \mathbf{r}\right) d%
\mathbf{r,}   \label{7}
\end{equation}
where $f_{\eta}$ is the distribution density of the vector%
\begin{equation*}
\eta\equiv\left\{ 
\begin{array}{ll}
\eta_{i}=\xi_{i}, & i\neq j, \\ 
\eta_{j}=\xi_{j}+\varepsilon\xi_{k}, & i=j,%
\end{array}
\right. 
\end{equation*}%
\begin{equation*}
f_{\eta}\left( \mathbf{r}\right) \equiv F\left( \varepsilon,\mathbf{r}%
\right) =\left\{ 
\begin{array}{ll}
f\left( r_{1},\ldots,r_{j-1},r_{j}-\varepsilon r_{k},\ldots,r_{N}\right) , & 
j\neq k, \\ 
\frac{1}{1+\varepsilon}f\left( r_{1},\ldots,r_{j-1},\frac{r_{j}}{%
1+\varepsilon},\ldots,r_{N}\right) , & j=k.%
\end{array}
\right. 
\end{equation*}

We see that%
\begin{equation*}
\frac{\partial F}{\partial\varepsilon}\left( \varepsilon,\mathbf{r}\right)
=\left\{ 
\begin{array}{ll}
-r_{k}\frac{\partial f}{\partial r_{j}}\left( r_{1},\ldots,r_{j}-\varepsilon
r_{k},\ldots,r_{N}\right) , & j\neq k, \\ 
\frac{-1}{\left( 1+\varepsilon\right) ^{2}}\frac{\partial}{\partial r_{j}}%
\left[ r_{j}f\left( r_{1},\ldots,\frac{r_{j}}{1+\varepsilon},\ldots
,r_{N}\right) \right] , & j=k.%
\end{array}
\right. 
\end{equation*}

Consider the function%
\begin{equation*}
\Psi\left( \varepsilon\right) \equiv\int\limits_{O\left( \mathbf{Z}\right) }%
\frac{\partial F}{\partial\varepsilon}\left( \varepsilon,\mathbf{r}\right) d%
\mathbf{r=-}\int\limits_{O_{\varepsilon}\left( \mathbf{Z}\right) }\Phi
_{kj}\left( \mathbf{r}\right) d\mathbf{r,}
\end{equation*}
where $O_{\varepsilon}\left( \mathbf{Z}\right) \equiv P_{\varepsilon
}O\left( \mathbf{Z}\right) $, $P_{\varepsilon}$ is a linear transformation
operator in $\mathbb{R}^{N}$: $P_{\varepsilon}\equiv\left\Vert
p_{il}\right\Vert $, $p_{il}=\delta_{il}-\alpha\left( \varepsilon\right)
\delta_{ij}\delta_{lk}$, 
\begin{equation*}
\left\{ 
\begin{array}{ll}
\alpha\left( \varepsilon\right) =\varepsilon, & k\neq j, \\ 
\alpha\left( \varepsilon\right) =\varepsilon/\left( 1+\varepsilon\right) , & 
k=j,%
\end{array}
\right. 
\end{equation*}
$\delta_{ij}$ is the Kronecker delta. We have%
\begin{equation*}
\left\vert \Psi\left( \varepsilon\right) \right\vert
\leq\int\limits_{R^{N}}\left\vert \Phi_{kj}\left( \mathbf{r}\right)
\right\vert d\mathbf{r,}
\end{equation*}
so the function $\Psi\left( \varepsilon\right) $ is integrable in a
neighborhood of the point $\varepsilon=0$. Let us prove that the function is
continuous at this point. Let denote for $R>0$%
\begin{equation*}
\begin{array}{ll}
K_{R}\equiv\left\{ \mathbf{r:}\left\vert \mathbf{r}\right\vert >R\right\} ,
& K_{R}^{\prime}\equiv\mathbb{R}^{N}\backslash K_{R}, \\ 
O_{\varepsilon,R}\left( \mathbf{Z}\right) \equiv O_{\varepsilon}\left( 
\mathbf{Z}\right) \cap K_{R}, & O_{\varepsilon,R}^{\prime}\left( \mathbf{Z}%
\right) \equiv O_{\varepsilon}\left( \mathbf{Z}\right) \cap K_{R}^{\prime},
\\ 
O_{0,R}\left( \mathbf{Z}\right) \equiv O\left( \mathbf{Z}\right) \cap K_{R},
& O_{0,R}^{\prime}\left( \mathbf{Z}\right) \equiv O\left( \mathbf{Z}\right)
\cap K_{R}^{\prime}.%
\end{array}
\end{equation*}
As the function $\Phi_{kj}\left( \mathbf{r}\right) $ is summable in $\mathbb{%
R}^{N}$ then for every $\delta>0$ there exist $R>0$ such that 
\begin{equation*}
\int\limits_{E}\left\vert \Phi_{kj}\left( \mathbf{r}\right) \right\vert d%
\mathbf{r<}\frac{\delta}{3}
\end{equation*}
for every measurable subset $E\subset K_{R}$. For $R>0$ do the estimation 
\begin{equation*}
\left\vert \Psi\left( \varepsilon\right) -\Psi\left( 0\right) \right\vert
=\left\vert \int\limits_{O_{\varepsilon}\left( \mathbf{Z}\right) }\Phi
_{kj}\left( \mathbf{r}\right) d\mathbf{r-}\int\limits_{O\left( \mathbf{Z}%
\right) }\Phi_{kj}\left( \mathbf{r}\right) d\mathbf{r}\right\vert \leq 
\end{equation*}%
\begin{equation*}
\leq\left\vert \int\limits_{O_{\varepsilon,R}\left( \mathbf{Z}\right)
}\Phi_{kj}\left( \mathbf{r}\right) d\mathbf{r}\right\vert \mathbf{+}%
\left\vert \int\limits_{O_{0,R}\left( \mathbf{Z}\right) }\Phi_{kj}\left( 
\mathbf{r}\right) d\mathbf{r}\right\vert \mathbf{+}\left\vert \int
\limits_{O_{\varepsilon,R}^{\prime}\left( \mathbf{Z}\right) }\Phi
_{kj}\left( \mathbf{r}\right) d\mathbf{r-}\int\limits_{O_{0,R}^{\prime
}\left( \mathbf{Z}\right) }\Phi_{kj}\left( \mathbf{r}\right) d\mathbf{r}%
\right\vert . 
\end{equation*}
The sum of the first two terms is less than $2\delta/3$. As the Lebesque
measure of the symmetric difference between bounded sets $O_{\varepsilon
,R}^{\prime}\left( \mathbf{Z}\right) $ and $O_{0,R}^{\prime}\left( \mathbf{Z}%
\right) $ tends to zero as $\varepsilon\rightarrow0$, then the last term is
less than $\delta/3$ for all sufficiently small $\varepsilon$. So we have $%
\left\vert \Psi\left( \varepsilon\right) -\Psi\left( 0\right) \right\vert
<\delta$ for all such $\varepsilon$. This proves the continuity of $%
\Psi\left( \varepsilon\right) $ at the point $\varepsilon=0$.

Listed above is sufficient for the possibility of differentiation of the
integral (\ref{7}) at the point $\varepsilon =0$ \cite[p. 132]{Sobol}: 
\begin{equation*}
\frac{\partial h_{k}}{\partial \mathbf{z}_{j}^{\ast }}=\int\limits_{O\left( 
\mathbf{Z}\right) }\frac{\partial F}{\partial \varepsilon }\left( 0,\mathbf{r%
}\right) ~d\mathbf{r=-}\int\limits_{O\left( \mathbf{Z}\right) }\Phi
_{kj}\left( \mathbf{r}\right) ~d\mathbf{r.}
\end{equation*}%
This proves the lemma.

For convex $\mathbf{Q}$ we have from (\ref{6})%
\begin{equation*}
\frac{\partial h_{k}}{\partial\mathbf{z}_{j}^{\ast}}=\int\limits_{\pi_{j}}%
\left( r_{k}f\right) \mid_{r_{j}=r_{j}^{\left( 2\right) }\left( \mathbf{p}%
_{j}\right) }^{r_{j}=r_{j}^{\left( 1\right) }\left( \mathbf{p}_{j}\right) }d%
\mathbf{p}_{j,}
\end{equation*}
where in the case of infinite values of $r_{j}^{\left( 1,2\right) }\left( 
\mathbf{p}_{j}\right) $ the substitution should be settled zero. Note that
using this representation of the derivative one can prove lemma 2 for more
wide assumptions about $f$ do not presuming the existence of the partial
derivatives.

\textbf{Theorem 2.} \textit{Suppose that }$N\geq n$\textit{, }$k\in 
\overline{1,N}$\textit{, set of vectors }$\left\{ \mathbf{z}_{1}^{\ast
},\ldots,\mathbf{z}_{n}^{\ast}\right\} \subset\mathbf{Z}$\textit{\ forms a
basis in }$R^{n}$\textit{, the density }$f$\textit{\ is differentiable with
respect to }$r_{j}$\textit{, and functions }$\Phi_{kj}\left( \mathbf{r}%
\right) \equiv\partial\left( r_{k}f\right) /\partial r_{j}$\textit{\ are
summable in }$R^{N}$\textit{. Then the function }$h_{k}\left( \mathbf{z}%
_{k}\right) $\textit{\ is differentiable at the point }$\mathbf{z}_{k}^{\ast
}$\textit{\ and its gradient is }%
\begin{equation}
\nabla h_{k}\left( \mathbf{z}_{k}^{\ast}\right)
=\sum\limits_{j=1}^{n}\sum\limits_{i=1}^{n}\frac{\partial h_{k}}{\partial%
\mathbf{z}_{j}^{\ast}}\left\langle \mathbf{e}^{j},\mathbf{e}%
^{i}\right\rangle \mathbf{z}_{i}^{\ast },   \label{8}
\end{equation}
\textit{where }$\left\{ \mathbf{e}^{1},\ldots,\mathbf{e}^{n}\right\} $%
\textit{\ is a dual basis for }$\left\{ \mathbf{z}_{1}^{\ast},\ldots ,%
\mathbf{z}_{n}^{\ast}\right\} $\textit{.}

\textbf{Prove.} If we prove the differentiability of the function, then
after decomposing of the gradient with basis $\left\{ \mathbf{e}^{i}\right\} 
$ and each of vectors $\mathbf{e}^{j}$ with basis $\left\{ \mathbf{z}%
_{i}^{\ast }\right\} $ \cite[p. 229]{Ilin} we receive%
\begin{equation*}
\nabla h_{k}=\sum\limits_{i=1}^{n}\left\langle \nabla h_{k},\mathbf{z}%
_{j}^{\ast }\right\rangle \mathbf{e}^{j}=\sum\limits_{j=1}^{n}\frac{\partial
h_{k}}{\partial \mathbf{z}_{j}^{\ast }}\sum\limits_{i=1}^{n}\left\langle 
\mathbf{e}^{j},\mathbf{e}^{i}\right\rangle \mathbf{z}_{i}^{\ast },
\end{equation*}%
that proves formula (\ref{8}).

To prove the differentiability of the function $h_{k}$ it is sufficiently to
check the continuity of the derivatives $\frac{\partial h_{k}}{\partial 
\mathbf{z}_{j}^{\ast}}\left( \mathbf{z}_{k}\right) $ at the point $\mathbf{z}%
_{k}^{\ast}$. After decomposition of an increment $\Delta \mathbf{z}_{k}$
with the basis $\left\{ \mathbf{z}_{i}^{\ast}\right\} $, $\Delta\mathbf{z}%
_{k}=\Delta_{1}\mathbf{z}_{1}^{\ast}+\ldots+\Delta _{n}\mathbf{z}_{n}^{\ast}$%
, $\Delta_{i}\equiv\left\langle \Delta\mathbf{z}_{k},\mathbf{e}%
^{i}\right\rangle $ we obtain by analogy with (\ref{7}) 
\begin{equation*}
\begin{array}{c}
\frac{\partial h_{k}}{\partial\mathbf{z}_{j}^{\ast}}\left( \mathbf{z}%
_{k}^{\ast}+\Delta\mathbf{z}_{k}\right) =-\int\limits_{O\left( \mathbf{Z}%
\right) }\frac{\partial}{\partial r_{j}}\left( r_{k}f\right) d\mathbf{r=} \\ 
=\frac{1}{1+\Delta_{k}}\int\limits_{O_{\Delta}\left( \mathbf{Z}\right)
}\Phi_{kj}\left( \mathbf{r}\right) d\mathbf{r=\Psi}\left( \Delta_{k}\right) ,%
\end{array}
\end{equation*}
where 
\begin{equation*}
\begin{array}{c}
f_{\Delta}\left( \mathbf{r}\right) =\frac{1}{1+\Delta_{k}}f\left( \mathbf{r}%
_{1}-\frac{\Delta_{1}}{1+\Delta_{k}}\mathbf{r}_{k},\ldots ,\mathbf{r}_{n}-%
\frac{\Delta_{n}}{1+\Delta_{k}}\mathbf{r}_{k},\mathbf{r}_{n+1},\ldots,%
\mathbf{r}_{N}\right) , \\ 
O_{\Delta}\left( \mathbf{Z}\right) \equiv P_{\Delta}O\left( \mathbf{Z}%
\right) ,\ \ \ \ \ P_{\Delta}\equiv\left\Vert p_{il}\right\Vert , \\ 
p_{il}\equiv\left\{ 
\begin{array}{ll}
\delta_{il}-\frac{\Delta_{1}}{1+\Delta_{k}}\delta_{lk}, & i=\overline{1,n}
\\ 
\delta_{il}, & i=\overline{n+1,N}%
\end{array}
\right. .%
\end{array}
\end{equation*}

We can to prove that $\Psi\left( \Delta_{k}\right) \rightarrow\Psi\left(
0\right) $\ as $\Delta\mathbf{z}_{k}\rightarrow0$ similar to lemma 2; this
proves the continuity of the derivatives and the theorem too.

It is clear that in the statement of the theorem one can change $\left\{ 
\mathbf{z}_{i}^{\ast}\right\} $ to any other basis $\mathbf{Z}$ in $\mathbb{R%
}^{n}$. If dimension of the linear span $\mathfrak{L}$ of $\mathbf{Z}$ is
less than $n$, then under formulated in the theorem properties of $f$ we can
guarantee only differentiability of $h_{k}$ over subspace $\mathfrak{L}$.

It is easy to construct the examples of functions $h_{k}\left( \mathbf{z}%
_{k}\right) $ which are not differentiable or even not continuous at the
given point $\mathbf{z}_{k}^{\ast}$ (see below). The reason of these may be
in absence of such properties of $\mathbf{Q}$ as bodility, convexity,
boundary smoothness.

\textbf{Example 1}. 
\begin{equation*}
\begin{array}{c}
n=3,\ N=2,\ \mathbf{Q}=\left\{ x_{1}=0,\left( x_{2}-2\right) ^{2}+\left(
x_{3}-2\right) ^{2}\leq1\right\} , \\ 
\mathbf{z}_{1}^{\ast}=\left( 0,0,1\right) ,\ \mathbf{z}_{2}^{\ast}=\left(
0,1,0\right) .%
\end{array}
\end{equation*}

\textbf{Example 2}. 
\begin{equation*}
\begin{array}{c}
n=2,\ N=1,\ \mathbf{Q}=\left\{ x_{1}\geq0,x_{2}\geq0,\right. \\ 
\left. 1\leq x_{1}^{2}+x_{2}^{2}\leq2\right\} ,\ \mathbf{z}^{\ast}=\left(
0,1\right) .%
\end{array}
\end{equation*}

\textbf{Example 3}. 
\begin{equation*}
\begin{array}{c}
n=2,\ N=1,\ \mathbf{Q}=\left\{ \left( x_{1}-7\right)
^{2}+x_{2}^{2}\leq25,\right. \\ 
\left. \left( x_{2}-7\right) ^{2}+x_{1}^{2}\leq25\right\} ,\ \mathbf{z}%
^{\ast}=\left( 1,1\right) .%
\end{array}
\end{equation*}

\textbf{Theorem 3.} \textit{Let function }$f$\textit{\ be smooth in }$R^{N}$%
\textit{, }$\mathbf{Q}$\textit{\ be a ball in }$R^{n}$\textit{. For every
set of vectors }$\mathbf{Z}$\textit{, which does not belong to hyperplane
holding the origin and touching the boundary of }$\mathbf{Q}$\textit{, the
functions }$h_{k}\left( \mathbf{z}_{k}\right) $,$\ k=\overline{1,N}$\textit{%
\ are differentiable at points }$\mathbf{z}_{k}^{\ast}\in\mathbf{Z}$\textit{.%
}

Using the analytical form of $\mathbf{Q}$, we can prove theorem 3 by direct
calculations. Note that for any set of vectors belonging to tangential plane
the functional equals zero; it can be maximal only for a degeneration
problem (see sec. \ref{nec}). Hence the statements of the theorem are
natural.

\section{Conditions of optimality}

\label{nec}

Denote by $\mathbf{Z}\left( \mathbf{u}\left( \cdot \right) \right) \equiv
\left\{ \mathbf{z}_{u}^{\left( i\right) },i=\overline{1,N}\right\} $ the set
of vectors in $\mathbb{R}^{n}$, described by formulas (\ref{zu}) for a given
control $\mathbf{u}\left( \cdot \right) \in \mathfrak{U}$ . The problem (\ref%
{1}) is called a) \textit{degenerate} if the functional $\varphi \left(
\cdot \right) $ is constant in $\mathfrak{U}$ ; b) \textit{regular for the
control} $\mathbf{u}^{\ast }\left( \cdot \right) \in \mathfrak{U}$ , if for
any $k=\overline{1,N}$ the functions $h_{k}\left( \mathbf{z}_{k}\right) $,
defined for the set $\mathbf{Z}\left( \mathbf{u}^{\ast }\left( \cdot \right)
\right) $, are differentiable at the corresponding points $\mathbf{z}%
_{u^{\ast }}^{\left( k\right) }$ (see sec. \ref{dif}).

\textbf{Lemma 3.} \textit{If }$\hat{x}_{0}\in Q$\textit{, then the
sufficient condition for }$u^{0}\left( \cdot\right) \in U$\textit{\ to be
optimal is that the set }$Z\left( \mathbf{u}^{0}\left( \cdot\right) \right) $%
\textit{\ be trivial (i.e. consists of zero vectors only). On the contrary,
if the problem is not degenerate and }$\hat{x}_{0}\not \in Q$\textit{, then
the necessary condition for }$u^{0}\left( \cdot\right) \in U$\textit{\ to be
optimal is that the set }$Z\left( \mathbf{u}^{0}\left( \cdot\right) \right) $%
\textit{\ be nontrivial.}

\textbf{\ Prove.} If $\mathbf{z}_{u^{0}}^{\left( i\right) }=0,i=\overline {%
1,N}$, then from (\ref{5}) we have $\varphi\left( \mathbf{u}^{0}\left(
\cdot\right) \right) =P\left( \mathbf{\hat{x}}_{0}\in\mathbf{Q}\right) $.
The corresponding value of the probability equals $1=\sup\varphi\left( 
\mathbf{u}\left( \cdot\right) \right) $ for $\mathbf{\hat{x}}_{0}\in\mathbf{Q%
}$ and $0=\inf\varphi\left( \mathbf{u}\left( \cdot\right) \right)
<\sup\varphi\left( \mathbf{u}\left( \cdot\right) \right) $ for $\mathbf{\hat{%
x}}_{0}\not \in \mathbf{Q}$. This proves both statements of lemma.

\textbf{Theorem 4.}\textit{\ Let }$\mathbf{u}^{0}\left( \cdot\right) $%
\textit{\ be optimal control, and the problem (\ref{1})\ is regular for the
control }$\mathbf{u}^{0}\left( \cdot\right) $\textit{. Then for almost all }$%
t\in\left[ 0,1\right] $%
\begin{equation}
\left\langle \mathbf{\theta}\left( t\right) ,B\left( t\right) \mathbf{u}%
^{0}\left( t\right) \right\rangle ={\max_{\mathbf{v\in V}}}\left\langle 
\mathbf{\theta}\left( t\right) ,B\left( t\right) \mathbf{v}\right\rangle , 
\label{9}
\end{equation}
\textit{where the conjugate function}%
\begin{equation*}
\mathbf{\theta}\left( t\right) =\sum\limits_{k=1}^{N}\chi_{k}\left( t\right) 
\mathbf{\theta}_{k}\left( t\right) , 
\end{equation*}%
\begin{equation*}
\left\{ 
\begin{array}{l}
\frac{d}{dt}\mathbf{\theta}_{k}=-A^{\prime}\left( t\right) \mathbf{\theta }%
_{k} \\ 
\mathbf{\theta}_{k}\left( 1\right) =\nabla h_{k}\left( \mathbf{z}%
_{u^{0}}^{\left( k\right) }\right) ,%
\end{array}
\right. 
\end{equation*}
and $\chi_{k}\left( t\right) $\textit{\ is the indicator of the set }$\left[
t_{k-1},t_{k}\right) ,\ k=\overline{1,N}$\textit{.}

\textbf{Prove.} Let $\tau\in\left( t_{k-1},t_{k}\right) $, $0<\varepsilon
<t_{k}-\tau$ 
\begin{equation*}
\mathbf{u}^{\varepsilon}\left( t\right) =\left\{ 
\begin{array}{ll}
\mathbf{v\in V}, & t\in\left( \tau,\tau+\varepsilon\right) , \\ 
\mathbf{u}^{0}\left( t\right) , & t\not \in \left( \tau,\tau+\varepsilon
\right) .%
\end{array}
\right. 
\end{equation*}

Taking into account that for such kind of variation only vector $\mathbf{z}%
_{u^{0}}^{\left( k\right) }$ varies, we have 
\begin{equation*}
\begin{array}{c}
\varphi\left( \mathbf{u}^{\varepsilon}\left( \cdot\right) \right)
-\varphi\left( \mathbf{u}^{0}\left( \cdot\right) \right) =h_{k}\left( 
\mathbf{z}_{u^{\varepsilon}}^{\left( k\right) }\right) -h_{k}\left( \mathbf{z%
}_{u^{0}}^{\left( k\right) }\right) = \\ 
=\left\langle \nabla h_{k}\left( \mathbf{z}_{u^{0}}^{\left( k\right)
}\right) ,\Delta_{\varepsilon}\mathbf{z}^{\left( k\right) }\right\rangle
+o\left( \left\vert \Delta_{\varepsilon}\mathbf{z}^{\left( k\right)
}\right\vert \right) ,%
\end{array}
\end{equation*}%
\begin{equation*}
\begin{array}{c}
\Delta_{\varepsilon}\mathbf{z}^{\left( k\right) }\equiv\mathbf{z}%
_{u^{\varepsilon}}^{\left( k\right) }-\mathbf{z}_{u^{0}}^{\left( k\right) }=
\\ 
=\Phi\left( 1\right) \int\limits_{\tau}^{\tau+\varepsilon}\Phi^{-1}\left(
s\right) B\left( s\right) \left( \mathbf{v-u}^{0}\left( s\right) \right) ds.%
\end{array}
\end{equation*}
Using the Lebesque theorem about differentiation of integrals we obtain for
almost all $\tau\in\left( t_{k-1},t_{k}\right) $ 
\begin{equation*}
\begin{array}{c}
\delta\varphi\equiv\underset{\varepsilon\rightarrow0}{\lim}\left(
\varphi\left( \mathbf{u}^{\varepsilon}\left( \cdot\right) \right)
-\varphi\left( \mathbf{u}^{0}\left( \cdot\right) \right) \right)
/\varepsilon= \\ 
=\left\langle \nabla h_{k}\left( \mathbf{z}_{u^{0}}^{\left( k\right)
}\right) ,\Phi\left( 1\right) \Phi^{-1}\left( \tau\right) B\left(
\tau\right) \left( \mathbf{v-u}^{0}\left( \tau\right) \right) \right\rangle =
\\ 
=\left\langle \mathbf{\theta}\left( \tau\right) ,B\left( \tau\right) \left( 
\mathbf{v-u}^{0}\left( \tau\right) \right) \right\rangle .%
\end{array}
\end{equation*}
Because the control $\mathbf{u}^{0}\left( \cdot\right) $ is optimal, the
variation $\delta\varphi\leq0$; this implies (\ref{9}). The theorem is
proved.


\begin{thebibliography}{9}
\bibitem{Smirn1} I.P. Smirnov, \textit{Control of the probability of the
system entry to a given region (in Russian).} Differential equations. 1990,
v. \textbf{26}, N 10, p. 1753-1758.

\bibitem{Smirn2} I.P. Smirnov, \textit{Necessary conditions for optimality
in the problem of entrance of stochastic control process to the given region
(in Russian).} Differential equations. 1990, v. \textbf{26}, N 11, p.
1943-1949.

\bibitem{Smirn3} I.P. Smirnov and A.Yu. Zorin, \textit{Optimal control of
the linear system with the noise at control parameter.} Proc. of the
international conf. "Mathematical algorithms-1", Nizhny Novgorod, NNSU,
1995, p. 108-110.

\bibitem{Gihman} I.I. Gihman, A.V. Skorohod, \textit{Stochastic differential
equations (in Russian).} Naukova dumka: Kiev, 1968.

\bibitem{Neve} Jaques Neveu, \textit{Bases math\'{e}matiques du calcul des
probabilit\'{e}s.} Masson et cie, Paris, 1964.

\bibitem{Sobol} S.L. Sobolev, \textit{The equations of the mathematical
physics (in Russian).} Nauka: Moscow, 1966.

\bibitem{Ilin} V.A. Il'in, A.G. Pozdn'yak, \textit{Linear algebra (in
Russian).} Nauka: Moscow, 1984.
\end{thebibliography}
\end{document}